\documentclass[12pt]{article}
\usepackage{amsfonts}
\usepackage{amsmath}
\usepackage{amssymb}  
\renewcommand{\baselinestretch}{1}
\begin{document}
\begin{center}{\large\bf Odd multiperfect numbers}
\vskip 0.5 cm {\bf Shi-Chao Chen\quad\quad Hao Luo} \vskip.5cm
Institute of Contemporary Mathematics\\ Department of Mathematics
and Information Sciences,\\ Henan University, Kaifeng, 475001,
CHINA\\  E-mail address: schen@henu.edu.cn

\end{center}
\parskip = 0.2cm
\topmargin=-0.1cm \oddsidemargin=0cm \evensidemargin=0cm
\renewcommand\baselinestretch{1.0}
\vskip.5cm

\noindent{\bf Abstract}

A natural number $n$ is called {\it multiperfect} or
{\it$k$-perfect} for integer $k\ge2$ if $\sigma(n)=kn$, where
$\sigma(n)$ is the sum of the positive divisors of $n$. In this
paper, we establish the structure theorem of odd multiperfect
numbers analogous as Euler's theorem on odd perfect numbers. We
prove the divisibility of the Euler part of odd multiperfect numbers
and  characterize the forms of odd perfect numbers $n=\pi^\alpha
M^2$ such that $\pi\equiv\alpha(\text{ mod }8)$. We also present
some examples to show the nonexistence of odd perfect numbers as
applications.


\vskip.5cm

\noindent{\bf 1. Introduction}\vskip.5cm

Let $n$ be a positive integer and $\sigma(n)$ be the sum of the
positive divisors of the natural number $n$. We say $n$ is {\it
multiperfect} or {\it$k$-perfect} if $\sigma(n)=kn$. For example, 6
is a 2-perfect number and 120 is a 3-perfect number. One can see [4]
or [9] for a survey of multiperfect numbers.

If $n$ is an even 2-perfect number $n$, then the well known
Euclid-Euler theorem states that the only even 2-perfect number is
$2^{p-1}(2^p-1)$, where $2^p-1$ is Mersenne prime. The odd case is
different. It is not known that whether odd $k$-perfect numbers
exist for any $k\ge2$. But some properties of odd $k$-perfect
numbers have been investigated. For odd 2-perfect number $n$, Euler
has shown that $n$ has the form $n=\pi^\alpha M^2$, where $\pi$ is
prime, $(\pi,\alpha)=1$ and $\pi\equiv\alpha\equiv1(\text{ mod }4)$.
$\pi^\alpha$ is called the Euler factor. Recently, Brogghan and Zhou
extended Euler's theorem to 4-perfect numbers, and then to
$2^k$-perfect numbers [1]. One goal of this paper is to obtain a
slightly more precise description of the structure of odd
$k$-perfect numbers. It turns out that an odd $k$-perfect number $n$
has the form $n={\it\Pi}M^2$ for which $({\it\Pi},M)=1$ and all
exponents of prime factors of ${\it \Pi}$ are odd (see Theorem 2.2
for details). ${\it\Pi}$ is called the {\it Euler part} of the
multiperfect number $n$ which is analogous as the Euler factor of
2-perfect numbers.

It is interesting to study the  Euler part of $k$-perfect numbers.
It turns out that the properties of the Euler part of $k$-perfect
numbers can be used to prove the nonexistence of odd $k$-perfect
numbers. For instance, Starni [10] recently proved that if an odd
2-perfect number $n$ has the form $n=\pi^\alpha
3^{2\beta}Q^{2\beta}$ with $(3,Q)=1$, then
$3^{2\beta}|\sigma(\pi^\alpha)$. Using this result, he showed that
if $\pi\equiv1(\text{ mod }12), \alpha\equiv1,9(\text{ mod }12)$,
then there does not exist odd 2-perfect numbers
$n=\pi^\alpha3^{2\beta}M^2$. This result was extended to odd
$2^k$-perfect numbers by Brogghan and Zhou [1]. We will show that
the prime 3 can be replaced with the Fermat prime under suitable
 conditions (see Corollary 3.2).

An early result of Starni [11] on the odd 2-perfect numbers is that
there is no odd 2-perfect numbers decomposable into primes all of
the type $\equiv1(\text{ mod }4)$ if $n=\pi^\alpha M^2$ and
$\pi\not\equiv\alpha(\text{ mod }8)$. In this paper, we will
characterize the forms of odd 2-perfect numbers such that
$\pi\equiv\alpha(\text{ mod }8)$. As a consequence, we extend
Starni's results and show the nonexistence of some forms of odd
perfect numbers.

\vskip.5cm

\noindent{\bf 2. Structure of multiperfect numbers}\vskip.5cm
Recently, Broughan and Zhou [1] extended Euler's structure theorem
of odd 2-perfect numbers to odd $4$-perfect numbers, and then to odd
$2^k$-perfect numbers. From the proof of the structure of odd
4-perfect number, they observed the following beautiful fact, which
is of independent interest.

\noindent{\bf Theorem} (Broughan and Zhou [1, Theorem 2.2]). {\it
For all odd primes $p$, powers $j\ge 1$ and odd exponents $e>0$, we
have}
$$2^j\|\sigma(p^e)\Longleftrightarrow 2^{j+1}\|(p+1)(e+1).$$

This theorem was proved by discussing different cases and
constructing some polynomials. Here we will give a short proof by an
element argument. More general, we have the following \vskip.5cm

\noindent{\bf Theorem 2.1.}  {\it Let $p$ be a prime and $e$ be a
positive integer. Let $\nu_2(m)$ be the highest power of $2$
dividing integer $m$. Then we have} $$\nu_2(\sigma(p^e))
=\left\{\begin{array}{ll}\nu_2(p+1)+\nu_2\left(\frac{e+1}{2}\right),
\quad &if\ p>2 \ and \ e\equiv1(\text{ mod }2),\\
0,&otherwise.\end{array} \right.\eqno{(1)}$$
$$\nu_2(\sigma(p^e)-1)=\left\{\begin{array}{ll}
0, &if\ p>2 \ and \ e\equiv1(\text{ mod }2),\\
\nu_2(p+1)+\nu_2\left(\frac{e}{2}\right), &if\ p>2 \ and \
e\equiv0(\text{ mod }2),
\\1,&if\ p=2.\end{array}
\right.\eqno{(2)}$$

\noindent{\bf Proof.} It is obvious that
$\sigma(p^e)=1+p+\cdots+p^e\equiv1(\text{ mod }2)$ when $p=2$ or the
case $p$ is odd and $e$ is even. Now we assume that $p$ and $e$ are
both odd. We write $r=\nu_2\left(\frac{e+1}{2}\right)$ and
$\frac{e+1}{2}=2^rs$ for some odd integer $s$. Then
\begin{align*}\sigma(p^e)&=\frac{p^{e+1}-1}{p-1}\\
&=\frac{(p^2)^{\frac{e+1}{2}}-1}{p-1}\\
&=\frac{(p^2)^{2^rs}-1}{p-1}\\
&=(p+1)\frac{(p^{2s})^{2^r}-1}{p^2-1}\\
&=(p+1)((p^{2s})^{2^{r-1}}+1)((p^{2s})^{2^{r-2}}+1)\cdots(p^{2s}+1)
\frac{p^{2s}-1}{p^{2}-1}.
\end{align*}
Since $p$ is odd, we have
$$(p^{2s})^{2^{r-i}}+1=(p^{s2^{r-i}})^2+1\equiv 2\ (\text{mod } 4),\quad 1\le i \le r.$$
It follows that
$$\nu_2((p^{2s})^{2^{r-i}}+1)=1,\quad 1\le i \le r.$$Note that
$$\frac{p^{2s}-1}{p^{2}-1}=1+p^2+p^4+\cdots+(p^2)^{s-1}\equiv s\equiv 1(\text{ mod } 2).$$
We get $$\nu_2\left(\frac{p^{2s}-1}{p^{2}-1}\right)=0.$$ Therefore
\begin{align*}\nu_2(\sigma(p^e))&=\nu_2(p+1)+\sum_{i=1}^r\nu_2((p^{2s})^{2^{r-i}}+1)
+\nu_2\left(\frac{p^{2s}-1}{p^{2}-1}\right)\\
&=\nu_2(p+1)+r\\
&=\nu_2(p+1)+\nu_2\left(\frac{e+1}{2}\right).\end{align*} The
formula (2) follows from (1) and the fact
$$\nu_2(\sigma(p^e)-1)=\nu_2(p\sigma(p^{e-1}))=\nu_2(p)+\nu_2(\sigma(p^{e-1})).$$
This complete the proof of Theorem 1. \quad$\Box$ \vskip.5cm

As an application of Theorem 2.1, we will establish the explicit
structure theorem of $k$-perfect numbers for any integer
$k\ge2$.\vskip.5cm

\noindent{\bf Theorem 2.2.} {\it Let $n$ be odd and $k$-perfect with
$\nu_2(k)\ge1$ and $s$ be any integer satisfying $1\le s\le
\nu_2(k)$. Then $n$ has the shape $$n=p_1^{e_1}p_2^{e_2}\cdots
p_{s}^{e_s}M^2,\eqno{(3)}$$ where $M$ is a positive integer, the
$p_i$ are primes with $(p_i,M)=1$ and the $e_j$ are odd positive
integers. If $\nu_2(k)-s$ has a nonnegative partition
$$\nu_2(k)-s=a_1+a_2+\cdots+a_s+b_1+b_2+\cdots+b_s,\ a_i\ge0,\ b_j\ge0,\eqno{(4)}$$
then the primes $p_1, \cdots p_s$ satisfy
$$p_i\equiv 2^{a_i+1}-1(\text{ mod }2^{{a_i+2}})$$ and the exponents  $e_1, \cdots, e_s$
satisfy $$e_j\equiv 2^{b_j+1}-1(\mbox{ mod }2^{{b_j+2}}).$$}
\vskip.5cm

Before proving the theorem, we give the definition of the Euler part
of odd $k$-perfect numbers.

\noindent{\bf Definition 2.3} {\it The Euler part of an odd
$k$-perfect number $n$ with the shape (3) is denoted by} $${\it
\Pi}:=p_1^{e_1}p_2^{e_2}\cdots p_{s}^{e_s}.$$ \vskip.5cm

\noindent{\bf Remark.}

(1) Theorem 2.2 shows that if $n_1$ is $k_1$-perfect, $n_2$ is
$k_2$-perfect and $\nu_2(k_1)=\nu_2(k_2)$, then $n_1$ and $n_2$ have
the same shapes. Therefore we only consider $2^k$-perfect numbers in
many situations.

(2) Note that there are $\nu_2(k)$ shapes of an odd $k$-perfect
number $n$ since $s$ can take $\nu_2(k)$ values and each $s$ gives a
shape of $n$ as (3).\vskip.5cm

\noindent{\bf Proof of Theorem 2.2.} Let
$n=p_1^{\alpha_1}p_2^{\alpha_2}\cdots p_r^{\alpha_r}$ where the
$p_i$ are distinct odd primes. Then $n$ is $k$-perfect implies that
$$\sigma(n)=\sigma(p_1^{\alpha_1})\sigma(p_2^{\alpha_2})\cdots\sigma(p_r^{\alpha_r})=kn.$$
Since $\nu_2(k)\ge1$, it follows from (1) of Theorem 1 that some
$\alpha_i$ must be odd. Therefore we can write $n$ as
$$n=p_1^{e_1}p_2^{e_2}\cdots p_{s}^{e_s}\cdot q_1^{2f_1}q_2^{2f_2}\cdots q_t^{2f_t},$$
where  $e_1,\cdots, e_s$ are odd positive integers, $f_1,\cdots,f_t$
are positive integers, $p_1, \cdots, p_s, q_1,\cdots, q_t$ are odd
primes. By Theorem 2.1 we have
\begin{align*}\nu_2(k)&=\nu_2(kn)\\&=\nu_2(\sigma(n))\\&=\sum_{i=1}^s\nu_2(\sigma(p_i^{e_i}))
+\sum_{j=1}^t\nu_2(\sigma(q_j^{2f_j}))\\
&=\sum_{i=1}^s\left(\nu_2(p_i+1)+\nu_2\left(\frac{e_i+1}{2}\right)\right)\\
&=s+\sum_{i=1}^s\left(\nu_2\left(\frac{p_i+1}{2}\right)
+\nu_2\left(\frac{e_i+1}{2}\right)\right).\end{align*}Therefore
$$\nu_2(k)-s=\sum_{i=1}^s\nu_2\left(\frac{p_i+1}{2}\right)
+\sum_{i=1}^s\nu_2\left(\frac{e_i+1}{2}\right).$$ Given a
nonnegative partition of $\nu_2(k)-s$ such that
$$\nu_2(k)-s=a_1+a_2+\cdots+a_s+b_1+b_2+\cdots+b_s,\ a_i\ge0,\
b_j\ge0,$$ we take $\nu_2\left(\frac{p_i+1}{2}\right)=a_i$ and
$\nu_2\left(\frac{e_i+1}{2}\right)=b_i, 1\le i\le s.$¡¡ Note that
\begin{align*} \nu_2\left(\frac{c+1}{2}\right)=d &\Longrightarrow
\frac{c+1}{2}=2^d(2l+1)\ for\ some\ integer\
l\ge0\\
&\Longrightarrow c=2^{d+2}l+2^{d+1}-1\\
&\Longrightarrow c\equiv 2^{d+1}-1(\mbox{ mod }2^{d+2}).
\end{align*}Theorem 2.2 follows immediately. \quad$\Box$\vskip.5cm

We give two examples to recover the well known Euler's theorem on
2-perfect numbers and Broughan and Zhou's structure theorem on
4-perfect numbers [1, Theorem 2.1].\vskip.5cm

\noindent{\bf Example 2.1.} Let $n$ be an odd $2$-perfect number.
Then $\nu_2(2)=1$, $s=1$. $n$ has the unique form $n=\pi^\alpha M^2$
with $\pi$ prime, $\alpha$ odd and $(\pi,M)=1$. By (4), $a_1=b_1=0$.
So $\pi\equiv\alpha\equiv1(\text{ mod }4).$\vskip.5cm

\noindent{\bf Example 2.2.} Let $n$ be an odd 4-perfect number. Then
$\nu_2(4)=2$.  $s=1$ or 2.

If $s=1$, then $n=p^eM^2$. By (4), $a_1=1, b_1=0$ or $a_1=0, b_1=1$.
Therefore $p\equiv3(\text{ mod }8), e\equiv1(\text{ mod }4)$ or
$p\equiv1(\text{ mod }4), e\equiv3(\text{ mod }8)$.

If $s=2$, then $n=p_1^{e_1}p_2^{e_2}M^2$. By (4),
$a_1=a_2=b_1=b_2=0$. Therefore $p_1\equiv p_2\equiv e_1\equiv
e_2\equiv1(\text{ mod }4)$.\vskip.5cm

\noindent{\bf 3. Euler part of $2^k$-perfect numbers}\vskip.5cm

Based on a result of McDaniel [7], Starni [10] proved that if an odd
2-perfect number $n$ has the form $n=\pi^\alpha
3^{2\beta}Q^{2\beta}$ with $(3,Q)=1$, then
$3^{2\beta}|\sigma(\pi^\alpha)$. This result was generalized to odd
$2^k$-perfect numbers by Broughan and Zhou [2, Theorem 2.6].
Recalling that ${\it \Pi}$ is the Euler part of an odd $k$-perfect
number. In the following Theorem 3.1 we will prove the divisible
result of $\sigma({\it \Pi})$ for a prime power. As a corollary, we
extend the results mentioned above.\vskip.5cm

\noindent{\bf Theorem 3.1.} {\it Let $n={\it\Pi}
q^{2\beta}\prod_{i=1}^sp_i^{2\beta_i}$ be an odd  $2^k$-perfect
number, where $\it{\Pi}$ is the Euler part of $n$, $q$ and the $p_i$
are distinct odd primes. Then
$$q^{2\beta}|\sigma(\it\Pi)\Longleftrightarrow
\left\{\begin{array}{ll}2\beta_i+1\not\equiv0(\text{ mod
}\text{ord}_q(p_i)),\quad&\text{ if }\ p_i\not\equiv1(\text{ mod
}q),\\2\beta_i+1\not\equiv0(\text{ mod }q),\quad&\text{ if }\
p_i\equiv1(\text{ mod }q),
\end{array}\right.$$ where $\text{ord}_q(m)$ is the order of $m$ in
the multiplicative group $(\mathbb{Z}/q\mathbb{Z})^*$.}

\noindent{\bf Proof.} By the definition of $2^k$-perfect,
$$2^kn=2^k{\it\Pi} q^{2\beta}\prod_{i=1}^sp_i^{2\beta_i}
=\sigma({\it\Pi})\sigma(q^{2\beta})\sigma\left(\prod_{i=1}^sp_i^{2\beta_i}\right)=\sigma(n).$$Therefore
$$q^{2\beta}|\sigma({\it\Pi})\Longleftrightarrow (q,\prod_{i=1}^s\sigma
(p_i^{2\beta_i}))=1\Longleftrightarrow (q,\sigma
(p_i^{2\beta_i}))=1,\ i=1,\cdots,s.$$If $p_i\equiv1(\text{ mod }q)$,
then
$$\sigma(p_i^{2\beta_i})\equiv 2\beta_i+1(\text{ mod }q).$$Therefore
$$(q,\sigma (p_i^{2\beta_i}))=1\Longleftrightarrow 2\beta_i+1\not\equiv0(\text{ mod }q).$$
If $p_i\not\equiv1(\text{ mod }q)$, then
\begin{align*}(q,\sigma
(p_i^{2\beta_i}))=1&\Longleftrightarrow
\left(q,\frac{p_i^{2\beta_i+1}-1}{p_i-1}\right)=1 \\
&\Longleftrightarrow p_i^{2\beta_i+1}\not\equiv1(\text{ mod
}q)\\&\Longleftrightarrow 2\beta_i+1\not\equiv0(\text{ mod
}\text{ord}_q(p_i)).\end{align*} The theorem follows.
\quad$\Box$\vskip.5cm

\noindent{\bf Corollary 3.2.}  {\it Let $n={\it\Pi}
q^{2\beta}\prod_{i=1}^sp_i^{2\beta_i}$ be an odd $2^k$-perfect as in
Theorem 3.2. If  $q$ is a Fermat prime, that is $q=2^{2^t}+1$ for
some integer $t\ge0$, and
$\prod_{i=1}^s(2\beta_i+1)\not\equiv0(\text{ mod }q)$, then
$q^{2\beta}|\sigma(\it\Pi).$}

\noindent{\bf Proof.} Note that the Euler function
$\phi(q)=2^{2^t}$. Lagrange's theorem implies that
$\text{ord}_q(p_i)|\phi(q)$. This show that $\text{ord}_q(p_i)$ is a
power of 2. The corollary follows.\quad$\Box$\vskip.5cm

\noindent{\bf Remark.} If $n=\pi^\alpha
q^{2\beta}\prod_{i=1}^sp_i^{2\beta}$ be an odd 2-perfect number with
$\alpha\equiv\pi\equiv1(\text{ mod }4)$, then it is known that
$\beta\neq2$ [6], $\beta\neq3$ [5], $\beta\neq5,12,17,24,62$ [8],
and $\beta\neq6,8,11,14,18$ [2]. In [8] McDaniel and Hagis
conjecture that there does not exist such 2-perfect numbers for any
positive integer $\beta$. Corollary 3.2 can be used to prove this
conjecture in some special cases. For example, if $q=5$,
$(2\beta+1,5)=1$, $\pi\equiv1(\text{ mod }20), \alpha\equiv13(\text{
mod }20)$, then by Corollary 3.2, $5|\sigma(\pi^\alpha)$, but
$\sigma(\pi^\alpha)\equiv1+\alpha\equiv4(\text{ mod }5)$. Therefore
there does not exist such odd 2-perfect numbers. In particular,
$n=41^{13}5^{2\beta}\prod_{i=1}^sp_i^{2\beta}$ can not be an odd
2-perfect number for any $\beta\not\equiv2(\text{ mod }5)$ and odd
primes $p_i$.\vskip.5cm

Let $n=\pi^\alpha\prod_ip_i^{2\beta_i}$ be an odd 2-perfect number.
$\pi^\alpha$, with $\pi\equiv\alpha\equiv1(\text{ mod }4)$, is the
Euler's factor. In [11] Starni proved the following results:

(a)\ {\it$\pi\equiv\alpha(\mbox{ mod }8)$ if each prime
$p_i\equiv1(\text{ mod }4)$.}

(b)\ {\it $\sigma(\pi^\alpha)/2$ cannot be prime if each prime
$p_i\equiv3(\text{ mod }4)$.}

The conclusion (a) was proved based on a result of Ewell [3]. We
will extend (a) and (b) in the following Theorem 3.3 and 3.4
respectively independent of  Ewell's result. As a consequence, we
obtain some results on the nonexistence of odd 2-perfect
numbers.\vskip.5cm

\noindent{\bf Theorem 3.3.} {\it Let $n=\pi^\alpha M^2$ be an odd
2-perfect number, with $\pi$ prime, $(\pi,M)=1$ and
$\pi\equiv\alpha\equiv1(\text{ mod }4)$. Then
\begin{align*}&\sigma(M^2)\equiv1(\text{ mod
}4)\Longleftrightarrow\pi\equiv\alpha(\text{ mod
}8),\tag{5}\\&\sigma(M^2)\equiv3(\text{ mod
}4)\Longleftrightarrow\pi\equiv\alpha+4(\text{ mod
}8).\tag{6}\end{align*} In particular, if
$n=\pi^\alpha\prod_ip_i^{2\beta_i}$ with $p_i\equiv1(\text{ mod }4)$
or $n=\pi^\alpha\prod_jq_j^{2\gamma_j}$ with $q_j\equiv3(\text{ mod
}4)$, then $$\pi\equiv\alpha(\mbox{ mod }8).$$}

\noindent{\bf Proof.} Using the fact that $\sigma(\pi^\alpha
M^2)=2\pi^\alpha M^2$ and $M, \alpha$ are odd, we find that
$$\pi^\alpha\equiv\pi\equiv\frac{\sigma(\pi^\alpha)}{2}\sigma(M^2)(\text{ mod
}8).\eqno{(7)}$$ Note that $\pi\equiv\alpha\equiv1(\text{ mod }4)$
implies
 $\pi^4\equiv\alpha^4\equiv1(\text{ mod }16)$. It follows that
\begin{align*}\sigma(\pi^\alpha)&=1+\pi+\cdots+\pi^\alpha\\
&=(1+\pi+\pi^2+\pi^3)(1+\pi^4+\pi^8+\cdots+\pi^{\alpha-5})+\pi^{\alpha-1}(1+\pi)\\
&\equiv(1+\pi)(1+\pi^2)\frac{\alpha-1}{4}+1+\pi(\text{ mod
}16).\end{align*} Hence we have
$$\frac{\sigma(\pi^\alpha)}{2}\equiv(1+\pi^2)\frac{1+\pi}{2}\frac{\alpha-1}{4}+\frac{1+\pi}{2}
\equiv\frac{1+\pi}{2}\frac{1+\alpha}{2}(\text{ mod }8).$$ It follows
from (7) that
$$\pi\equiv\frac{\pi+1}{2}\frac{\alpha+1}{2}\sigma(M^2)(\text{ mod }8).$$
Since $\frac{1+\pi}{2}$ is odd, we have
$$\frac{\pi(\pi+1)}{2}\equiv\left(\frac{\pi+1}{2}\right)^2\frac{\alpha+1}{2}\sigma(M^2)
\equiv\frac{\alpha+1}{2}\sigma(M^2)(\text{ mod }8).\eqno{(8)}$$ If
$\sigma(M^2)\equiv1(\text{ mod }8)$, then
$\frac{\pi(\pi+1)}{2}\equiv\frac{\alpha+1}{2}(\text{ mod }8)$
implies that $$\pi(\pi+1)\equiv\alpha+1(\text{ mod }16).$$ Recall
that $\pi\equiv\alpha\equiv1(\text{ mod }4)$. It is easy to find the
solutions $(\pi, \alpha)(\text{ mod } 16)$ are
$$(1,1);\quad (5,13);\quad(9,9);\quad (13,5).$$ In particular, we
get $$\pi\equiv\alpha(\text{ mod }8).$$ Similarly, by (8), one can
find that
\begin{align*}\sigma(M^2)\equiv5(\text{ mod
}8)&\Longrightarrow\pi(\pi+1)\equiv5(\alpha+1)(\text{ mod }16)\\
&\Longrightarrow(\pi,\alpha)(\text{ mod }16)=(1,9);(5,5);(9,1);(13,13)\\
&\Longrightarrow\pi\equiv\alpha(\text{ mod }8).\end{align*}
\begin{align*}\sigma(M^2)\equiv3(\text{ mod
}8)&\Longrightarrow\pi(\pi+1)\equiv3(\alpha+1)(\text{ mod }16)\\
&\Longrightarrow(\pi,\alpha)(\text{ mod }16)=(1,5);(5,9);(9,13);(13,1)\\
&\Longrightarrow\pi\equiv\alpha+4(\text{ mod }8).\end{align*}
\begin{align*}\sigma(M^2)\equiv7(\text{ mod
}8)&\Longrightarrow\pi(\pi+1)\equiv7(\alpha+1)(\text{ mod }16)\\
&\Longrightarrow(\pi,\alpha)(\text{ mod }16)=(1,13);(5,1);(9,5);(13,9)\\
&\Longrightarrow\pi\equiv\alpha+4(\text{ mod }8).\end{align*} This
prove (5) and (6).

If we write $M^2=\prod_ip_i^{2\beta_i}\prod_j q_j^{2\gamma_j}$ where
$p_i\equiv1(\text{ mod }4), q_j\equiv3(\text{ mod } 4)$, then by (2)
of Theorem 1 $\sigma(q_j^{2\gamma_j})\equiv1(\text{ mod }4)$. Hence
$$\sigma(M^2)=\prod_i\sigma(p_i^{2\beta_i})\prod_j\sigma(q_j^{2\gamma_j})
\equiv\prod_i(2\beta_i+1)(\text{ mod } 4).\eqno{(9)}$$ Clearly, if
$\beta_i=0$ for all $i$, that is $M^2=\prod_j q_j^{2\gamma_j}$, then
$\sigma(M^2)\equiv1(\text{ mod } 4)$, and (4) implies that
$\pi\equiv\alpha(\text{ mod }8).$

If  $M^2=\prod_ip_i^{2\beta_i}$ with $p_i\equiv1(\text{ mod }4)$,
then
$$2n=2\pi^\alpha\prod_ip_i^{2\beta_i}=\sigma(\pi^\alpha)\sigma(M^2)=\sigma(n).$$
This shows that each prime factors of $\sigma(M^2)$ is
$\equiv1(\text{mod }4)$. Hence we have $\sigma(M^2)\equiv1(\text{mod
}4)$ and (5) implies that $\pi\equiv\alpha(\text{ mod
}8).$\quad$\Box$\vskip.5cm

\noindent{\bf Remark.} By (6), (9) and (2) of Theorem 2.1, it is
easy to see that $\pi\equiv\alpha+4(\text{ mod }8)$ if and only if
the number of prime factors $p^e$  of $M^2$ with $p^e\parallel M^2$,
$p\equiv1(\text{ mod } 4)$ and $e\equiv2(\text{ mod } 4)$ is
odd.\vskip.5cm

\noindent{\bf Theorem 3.4} {\it Let $n={\it\Pi}M^2$ be an odd
$2^k$-perfect number, where $\it{\Pi}$ is the Euler part of $n$. If
all prime factors of $M$ are $\equiv3(\text{ mod }4)$ and
${\it\Pi}=p_1^{e_1}\cdots p_s^{e_s}q_1^{f_1}\cdots q_{2t}^{f_{2t}}$
satisfies $(\sigma({\it\Pi}),p_1\cdots p_s)=1$, where the primes
$p_i\equiv1(\text{ mod }4)$, $q_j\equiv3{\text{ mod }4}$ and integer
$t\ge0$, then
$$\Omega\left(\frac{\sigma({\it\Pi})}{2^k}\right)\equiv0(\text{ mod }2),$$
where $\Omega(\sigma({\it\Pi})/2^k)$  is the total number of prime
factors of $\sigma({\it\Pi})/2^k$.}

\noindent{\bf Proof.} By the definition of $2^k$-perfect, we have
$$p_1^{e_1}\cdots
p_s^{e_s}q_1^{f_1}\cdots
q_{2t}^{f_{2t}}M^2=\frac{\sigma({\it\Pi})}{2^k}\sigma(M^2).$$ Since
$(\sigma({\it\Pi}),p_1\cdots p_s)=1$, we deduce that all prime
factors of $\frac{\sigma({\it\Pi})}{2^k}$ are $\equiv3(\text{ mod
}4)$. Note that $\sigma(M^2)\equiv1(\text{ mod }4)$. It follows that
$$1\equiv p_1^{e_1}\cdots
p_s^{e_s}q_1^{f_1}\cdots
q_{2t}^{f_{2t}}M^2=\frac{\sigma({\it\Pi})}{2^k}\sigma(M^2)
\equiv(-1)^{\Omega(\frac{\sigma({\it\Pi})}{2^k})}(\text{ mod }4).$$
Therefore
$$\Omega\left(\frac{\sigma({\it\Pi})}{2^k}\right)\equiv0(\text{ mod
}2).\quad\Box$$

If $n=\pi^\alpha M^2$ is an odd 2-perfect number such that all prime
factors of $M$ are $\equiv3(\text{ mod }4)$, then Theorem 3.5
implies that
$$\Omega\left(\frac{\sigma(\pi^\alpha)}{2}\right)\equiv0(\text{ mod
}2).$$ In particular, this show that $\frac{\sigma(\pi^\alpha)}{2}$
can not be a prime. We can use this fact to prove the nonexistence
of 2-perfect numbers. For example, $\pi=209$,
$\sigma(\pi)=210=2\cdot3\cdot5\cdot7$ and
$\Omega(\frac{\sigma(\pi)}{2})=3$. $\pi=30029$,
$\sigma(30029)=30030=2\cdot3\cdot5\cdot7\cdot11\cdot13$ and
$\Omega(\frac{\sigma(\pi)}{2})=5$. It follows that there does not
exist 2-perfect number $n=209M^2$ and $n=30029M^2$ with all prime
factors of $M$ are $\equiv3(\text{ mod }4)$.\vskip.5cm

\noindent{\bf Acknowledgements} The authors thank the Natural
Science Foundation of China (Grant No.11026080 ) and the Natural
Science Foundation of Education Department of Henan Province (Grant
No. 2009A110001).\vskip.5cm

\noindent{\bf References}\vskip.5cm

\begin{itemize}

\item[{[1]}] K. Broughan and Q. Zhou, Odd multiperfect numbers of abundancy 4, J.
Number Theory 128(2008) 1566-1575.

\item[{[2]}] G.L. Cohen and R.J.Williams, Extensions of some results concerning
odd perfect numbers, Fibonacci Quart. 23 (1985) 70-76.

\item[{[3]}] J.A.Ewell, On the multiplicative structure of odd
perfect numbers, J. Number Theory 12 (1980) 339-342.

\item[{[4]}] A.Flammenkamp, The multiply perfect numbers page,
http://wwwhomes.uni-bielefeld.de/achim/mpn.html.

\item[{[5]}] P.Hagis Jr., W.L. McDaniel, A new result concerning the structure of
odd perfect numbers, Proc. Amer. Math. Soc. 32 (1972) 13-15.

\item[{[6]}] H.-J. Kanold, Untersuchungen $\ddot{u}$ber ungenrade vollkommene
Zahlen, J. Reine Angew. Math. 183 (1941) 98-109.

\item[{[7]}] W.L. McDaniel, The nonexistence of odd perfect numbers of a certain
form, Arch. Math. 21(1970) 52-53.

\item[{[8]}] W.L. McDaniel, P.Hagis Jr., Some results concerning the nonexistence
of odd perfect numbers of the form $p^\alpha M^{2\beta}$,Fibonacci
Quart. 13 (1975) 25-28.

\item[{[9]}] R.Sorli, Multiperfect numbers,
http://www-staff.math.uts.edu.au/~rons/mpfn/mpfn.htm.

\item[{[10]}] P.Starni, On some properties of the Euler's factor of certain odd
perfect numbers, J. Number Theory 116 (2006) 483-486.

\item[{[11]}] P.Starni, On the Euler's factor of an odd perfect number, J. Number
Theory 37 (1991) 366-369.

\end{itemize}

\end{document}